\documentclass{article}
\usepackage[pdftex]{}
\usepackage{graphicx}
\usepackage{stmaryrd}
\usepackage{bm}
\usepackage{amsmath}
\usepackage{amsfonts}
\usepackage{amsthm}
\usepackage{mathrsfs}
\usepackage{multirow,bigdelim}
\usepackage{comment}
\usepackage[english]{babel}
\usepackage{url}

\theoremstyle{definition}
\newtheorem{theo}{Theorem}[section]
\newtheorem*{theo*}{Theorem}

\newtheorem*{defi*}{Definition}
\newtheorem{lemm}[theo]{Lemma}
\newtheorem*{fact*}{Fact}
\newtheorem*{ex*}{Example}
\newtheorem*{re*}{Remark}
\newtheorem{prop}[theo]{Proposition}
\newtheorem*{ack*}{Acknowledgement}
\newtheorem*{A*}{Theorem A}
\newtheorem*{B*}{Theorem B}

\makeatletter
\newcommand{\xRightarrow}[2][]{%
\ext@arrow 0055{\Rightarrowfill@}{#1}{#2}%
}
\def\Rightarrowfill@{\arrowfill@\Relbar\Relbar\Rightarrow}
\newcommand{\xLeftarrow}[2][]{%
\ext@arrow 0055{\Leftarrowfill@}{#1}{#2}%
}
\def\Leftarrowfill@{\arrowfill@\Leftarrow\Relbar\Relbar}
\newcommand{\xLongleftrightarrow}[2][]{%
\ext@arrow 0055{\llrafill@}{#1}{#2}%
}
\def\llrafill@{\arrowfill@\Leftarrow\Relbar\Rightarrow}
\makeatother

\makeatletter
\newenvironment{Proof}[1][\Proofname]{\par
  \normalfont
  \topsep6\p@\@plus6\p@ \trivlist
  \item[\hskip\labelsep{\bfseries #1}\@addpunct{\bfseries.}]\ignorespaces
}{%
  \endtrivlist
}

\newcommand{\Proofname}{Proof}
\makeatother
\makeatletter
\def\BOXSYMBOL{\RIfM@\bgroup\else$\bgroup\aftergroup$\fi
  \vcenter{\hrule\hbox{\vrule height.85em\kern.6em\vrule}\hrule}\egroup}
\makeatother
\newcommand{\BOX}{%
  \ifmmode\else\leavevmode\unskip\penalty9999\hbox{}\nobreak\hfill\fi
  \quad\hbox{\BOXSYMBOL}}

\newcommand\QED{\BOX}

\def\Rv{R_v}
\def\DM{\Delta_{\mathcal{M}}}
\def\Wvs{W_{v,s}}
\def\Wvone{W_{v,1}}
\def\Wvtwo{W_{v,2}}
\def\Vvs{V_{v,s}}
\def\Vvone{V_{v,1}}

\pdfoutput=1

\begin{document}

\title{The smooth torus orbit closures in the Grassmannians}
\author{Masashi Noji and Kazuaki Ogiwara}
\date{}
\maketitle

\begin{abstract}
It is known that for the natural algebraic torus actions on the Grassmannians, the closures of torus orbits are toric varieties, and that these toric varieties are smooth if and only if the corresponding matroid polytopes are simple. We prove that simple matroid polytopes are products of simplices and smooth torus orbit closures in the Grassmannians are products of complex projective spaces. Moreover, it turns out that the smooth torus orbit closures are uniquely determined by the corresponding simple matroid polytopes. 
\end{abstract}

\begin{center}
{\Large \textbf{\sf Introduction}}
\end{center}

\fontsize{10pt}{0.55cm}\selectfont
Let $G_k(\mathbb{C}^n)$ be the Grassmannian consisting of all $k$-dimensional linear subspaces of the $n$-dimensional complex linear space $\mathbb{C}^n$ and set $T = (\mathbb{C^*})^n$. The coordinatewise multiplication of $T$ on $\mathbb{C}^n$ induces a $T$-action on $G_k(\mathbb{C}^n)$. 
The study of $T$-orbits in $G_k(\mathbb{C}^n)$ was initiated by Gel'fand \cite{g}, Gel'fand-Serganova \cite{gs} and Aomoto \cite{a} in connection with hypergeometric functions and Gel'fand-MacPherson \cite{gm} constructed the correspondence between polytopes and $T$-orbit closures to get differential forms on the (real) Grassmannians. This correspondence was developed by Gel'fand-Goresky-MacPherson-Serganova in \cite{ggms}. They discovered a connection between $T$-orbits, the theory of polytopes, and the theory of combinatorial geometries. Recently, Buchstaber-Terzi\'c \cite{bs1}, \cite{bs2} investigated $T$-orbits and their closures in detail especially when $(n,k)=(4,2), (5,2)$. 

It is known that the closure of any $T$-orbit in $G_k(\mathbb{C}^n)$ is normal and hence a toric variety. 
These $T$-orbit closures are not necessarily smooth. In this paper we will completely identify \emph{smooth} $T$-orbit closures in $G_k(\mathbb{C}^n)$ for an arbitrary value of $(n,k)$. 

If $(n_1,\dots,n_r)$ is a partition of $n$ and $(k_1 ,\dots, k_r)$ is a partition of $k$, then the product of Grassmannians 
\begin{equation} \label{eq:0-1}
G_{k_1}(\mathbb{C}^{n_1})\times G_{k_2}(\mathbb{C}^{n_2})\times \dots\times G_{k_r}(\mathbb{C}^{n_r})
\end{equation}
is naturally sitting in $G_k(\mathbb{C}^n)$ and stable under the $T$-action. As is easily observed, the product above is a $T$-orbit closure if and only if each $G_{k_i}(\mathbb{C}^{n_i})$ is a point or a complex projective space, i.e. $k_i\in \{0,1,n_i-1,n_i\}$ for each $i$. This provides examples of smooth $T$-orbit closures in $G_k(\mathbb{C}^n)$. Our main result in this paper stated below says that these are essentially all smooth $T$-orbit closures in $G_k(\mathbb{C}^n)$. 

\begin{A*}[see Theorem 2.1] \label{theo:0-1}
If a $T$-orbit closure in the Grassmannian $G_k(\mathbb{C}^n)$ is smooth, then it is of the form \eqref{eq:0-1} with $k_i\in \{0,1,n_i-1,n_i\}$ for each $i$ up to permutations of coordinates of $\mathbb{C}^n$ and hence isomorphic to a product of complex projective spaces. 
\end{A*}

The proof of Theorem A reduces to combinatorics on polytopes, which we shall explain. 
Gel'fand- 

\vspace{0.1in}
\hrule width 3cm
\vspace{0.05in}

\fontsize{8pt}{0.4cm}\selectfont \noindent 
{\sl Date}: December 10, 2018 



\noindent
Ogiwara and Noji were partially supported by the bilateral program ``Topology and geometry of torus actions, cohomological rigidity, and hyperbolic manifolds'' between JSPS and RFBR.

\fontsize{10pt}{0.55cm}\selectfont
\noindent Goresky-MacPherson-Serganova \cite{ggms} introduced a map $\mu\colon G_k(\mathbb{C}^n) \to \mathbb{R}^n$ called a moment map. It is invariant under the $T$-action on $G_k(\mathbb{C}^n)$ and the image $\mu(G_k(\mathbb{C}^n))$ is the hypersimplex $\Delta_{n,k}$ that is the convex hull of $\binom{n}{k}$ points in $\mathbb{R}^n$ with $k$ $1$'s and $n-k$ $0$'s in the coordinates. 
They show that the image of a $T$-orbit closure $\overline{\mathcal{O}}$ in $G_k(\mathbb{C}^n)$ by $\mu$ is a polytope sharing vertices and edges with $\Delta_{n,k}$. Such a polytope is called an $(n,k)$-hypersimplex in \cite{gs} but called a matroid polytope nowadays.  
A key fact we will use is that the $T$-orbit closure $\overline{\mathcal{O}}$ is smooth if and only if the matroid polytope $\mu(\overline{\mathcal{O}})$ is simple (cf. [6, Proposition 1 in p.150]). 

We note that if $(n_1,\dots,n_r)$ is a partition of $n$ and $(k_1 ,\dots, k_r)$ is a partition of $k$, then the product of hypersimplices 
\begin{equation} \label{eq:0-2}
\Delta_{n_1,k_1}\times \Delta_{n_2,k_2}\times \dots \times \Delta_{n_r,k_r}
\end{equation}
is naturally sitting in $\Delta_{n,k}$ and it is simple if and only if each $\Delta_{n_i,k_i}$ is a point or a simplex, i.e. $k_i\in \{0,1,n_i-1,n_i\}$ for each $i$. It turns out that these are essentially all simple matroid polytopes in $\Delta_{n,k}$ as stated below and this combinatorial result implies Theorem A above. 

\begin{B*}[see Theorem 1.3] \label{theo:0-2} 
If a matroid polytope in the hypersimplex $\Delta_{n,k}$ is simple, then it is of the form \eqref{eq:0-2} with $k_i\in \{0,1,n_i-1,n_i\}$ for each $i$ up to permutations of coordinates of $\mathbb{R}^n$ and hence isomorphic to a product of simplices. 
\end{B*}

This paper is organized as follows. The first section is devoted to the proof of Theorem B above. To prove it, we associate a graph to each vertex $v$ of a matroid polytope by using the edge vectors emanating from the vertex $v$. This graph is an analogue of that introduced in \cite{em}. The key fact we use is that the $v$ is a simple vertex if and only if the graph has no cycle. 
In the second section, we prove Theorem A by applying Theorem B through the moment maps. 

\begin{ack*}
We thank Professor Mikiya Masuda for bringing the papers \cite{ggms}, \cite{gs}, \cite{bs1}, \cite{bs2} to our attention and for his support throughout the project.
\end{ack*}

\section{Simple matroid polytopes}

For a set $S$, let $\mathbb{R}^S$ be the set of all real-valued functions on $S$ and let $[n]$ be the set of all natural numbers from $1$ to $n$. When $S = [n]$, $\mathbb{R}^{[n]}$ can be identified with $\mathbb{R}^n$. For $B \subset S$, let $\delta_B$ be the real-valued function on $S$ defined by
\begin{equation}
(\delta_B)(i )= 
\begin{cases}
1 \quad (i \in B) \\ \notag
0 \quad (i \notin B). \notag
\end{cases}
\end{equation}

For an integer $k$ between $0$ and $n$, let $\Delta_{S,k}$ be the convex hull of $\{\delta_B | \ B\subset S,|B| = k\}$ where the cardinality of the empty set $\emptyset$ is understood to be $0$. 
When $S = [n]$, this convex polytope is called a \textbf{hypersimplex} and denoted by $\Delta_{n,k}$. 
The hypersimplex $\Delta_{n,k}$ is a point when $k=0$ or $n$ and is of dimension $n-1$ when $1\le k\le n-1$. The vertices of $\Delta_{n,k}$ are $\delta_B$'s with $|B|=k$ and two vertices $\delta_B$ and $\delta_{B'}$ of $\Delta_{n,k}$ are joined by an edge if and only if $|B\cap B'|=k-1$. Therefore, $\Delta_{n,k}$ has $\binom{n}{k}$ vertices and $k(n-k)$ edges at each vertex. In particular, $\Delta_{n,k}$ is not simple for $2\le k\le n-2$. Note that when $k=1$ or $n-1$, $\Delta_{n,k}$ is a simplex and hence simple. 

\newpage
\begin{ex*}
The hypersimplex $\Delta_{4,2}$ is the regular octahedron:
\begin{center}
\includegraphics[width = 3 cm]{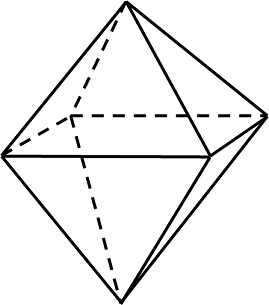}

Figure 1.1\quad The hypersimplex $\Delta_{4,2}$.
\end{center}
\end{ex*}

\begin{defi*}[Matroid and matroid polytope]
A nonempty collection $\mathcal{M}$ of $k$-element subsets of $S$ is called a matroid of rank $k$ on $S$ if it satisfies the following: if $I$ and $J$ are distinct members of $\mathcal{M}$ and $i \in I\setminus J$, then there exists an element $j \in J \setminus I$ such that $(I \setminus \{i\}) \cup \{j\} \in \mathcal{M}$. 
A matroid polytope (or a matroid basis polytope) is the convex hull of $\delta_I$'s $(I\in \mathcal{M})$ in $\mathbb{R}^S$ for a matroid $\mathcal{M}$ and denoted by $\DM$. 
\end{defi*}

\begin{re*}
Any matroid polytope in $\mathbb{R}^S$ has dimension at most $|S| - 1$ because the points $\delta_I$ with $|I|=k$ for a fixed $k$ lie on a hyperplane of $\mathbb{R}^S$. 
\end{re*}

We recall two facts on matroid polytopes used in this paper. 
The first one is the following. For a convex polytope $P$, we denote by $V(P)$ (resp. $E(P)$) the set of all vertices (resp. edges) of $P$. 

\begin{prop}[\cite{ggms}, 4.1.THEOREM] \label{prop:1-1}
Let $\DM$ be the matroid polytope of a matroid $\mathcal{M}$ of rank $k$ on $S$. Then 
\[
\begin{split}
V(\DM)&=\{\delta_I \mid I\in \mathcal{M}\},\\
E(\DM)&=\{ \delta_I\delta_J \mid I,J\in\mathcal{M},\ |I\cap J|=k-1\}\\
&=\{ \delta_I\delta_J \mid I,J\in\mathcal{M},\ \delta_I\delta_J\in E(\Delta_{S,k})\}
\end{split}
\]
where $\delta_I\delta_J$ denotes the edge between the vertices $\delta_I$ and $\delta_J$. 
Moreover, a convex polytope $\Delta \subset \Delta_{S,k}$ is a matroid polytope if and only if 
$V(\Delta) \subset V(\Delta_{S,k})$ and $E(\Delta) \subset E(\Delta_{S,k})$. 
\end{prop}

\begin{ex*}
Let $S = [4]$ and $k = 2$. Then the polytope on the left side of Figure 1.2 is not a matroid polytope because the polytope has an edge which is not an edge of $\Delta_{S,k}$. On the other hand, the polytope on the right side of Figure 1.2 is a matroid polytope.
\begin{center}
\includegraphics[width = 6 cm]{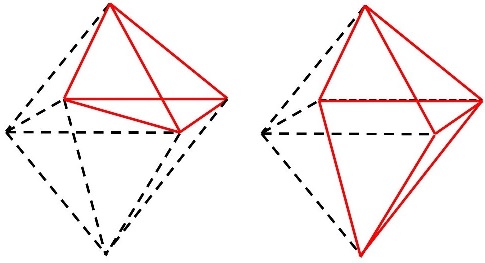}

Figure 1.2\quad Not a matroid polytope and a matroid polytope 
\end{center}
\end{ex*}

The second fact on matroid polytopes used in this paper is the following, which can be found in the proof of \cite[Proposition 4]{gs} for example. 

\begin{prop}[\cite{gs}, Proposition 4] \label{prop:1-2}
If $\mathcal{M}$ is a matroid of rank $k$ on $S$ and $\dim\DM=|S|-q$, then there are a partition $(J_1,J_2,\dots,J_q)$ of $S$ and a matroid $\mathcal{M}_i$ of rank $k_i$ on $J_i$ for each $i=1,2,\dots,q$ such that $(k_1,k_2,\dots,k_q)$ is a partition of $k$ and 
\begin{equation} \label{eq:product}
\DM=\Delta_{\mathcal{M}_1}\times \Delta_{\mathcal{M}_2}\times\dots \times \Delta_{\mathcal{M}_q} \subset \Delta_{S,k}\subset \mathbb{R}^{J_1} \times \mathbb{R}^{J_2} \times \cdots \times \mathbb{R}^{J_q} = \mathbb{R}^S
\end{equation}
where $\dim \Delta_{\mathcal{M}_i}=|J_i|-1$ for each $i=1,2,\dots,q$.
\end{prop}

Our main result in this section is the following:

\begin{theo}\label{theo1}
Let the situation be as in Proposition~\ref{prop:1-2}. If $\DM$ is simple, then $\Delta_{\mathcal{M}_i}=\Delta_{J_i,k_i}$ and $k_i \in \{0,1,|J_i|-1,|J_i|\}$ for each $i=1,2,\dots,q$ in \eqref{eq:product}. 
In particular, $\DM$ is a product of simplices if it is simple. 
\end{theo}

\begin{re*}
Any $2$-face of a matroid polytope is combinatorially equivalent to a triangle or a square because any face of a matroid polytope is a matroid polytope and any two dimensional matroid polytope is combinatorially equivalent to a triangle or a square. On the other hand, it is known that if any $2$-face of a simple polytope is combinatorially equivalent to a triangle or a square, then the simple polytope is combinatorially equivalent to a product of simplices (see \cite[Appendix]{ym} for example). It follows from this argument that a simple matroid polytope is combinatorially equivalent to a product of simplices. However, the statement in Theorem~\ref{theo1} is stronger than this combinatorial statement because the theorem describes the simple polytope $\DM$ precisely and also shows how it lies in $\Delta_{S,k}$. 
\end{re*}

The rest of this section is devoted to the proof of Theorem~\ref{theo1}. We may take $S=[n]$. Note that the product $P_1\times P_2\times \cdots \times P_q$ of convex polytopes $P_i$ is simple if and only if every $P_i$ is simple. Therefore, it suffices to prove the theorem when $\DM$ is of maximal dimension, i.e. $\dim \DM=n-1$ by Proposition~\ref{prop:1-2}.

Let $\Phi = \{\delta_i-\delta_j \mid i \neq j \in [n] \}$. For $R \subset \Phi$ with $R\cap (-R)=\emptyset$, we define a graph $\Gamma(R)$ as follows: the vertex set is $[n]$ and
\begin{center}
$\{i,j\}$ is an edge of $\Gamma(R)$ $\Leftrightarrow \delta_i - \delta_j \in R$ or $\delta_j - \delta_i \in R$. 
\end{center}
Note that the graph $\Gamma(R)$ is simple, i.e. it has no multiple edges and loops.

\begin{ex*}
If $n = 4, R_1 = \{\delta_1-\delta_3, \delta_1-\delta_4,\delta_2-\delta_3,\delta_2-\delta_4\}, R_2 = \{\delta_1 - \delta_3, \delta_2 - \delta_3\}$, then $\Gamma(R_1), \Gamma(R_2)$ are the following graphs:
\begin{center}
\includegraphics[width = 2.7 cm]{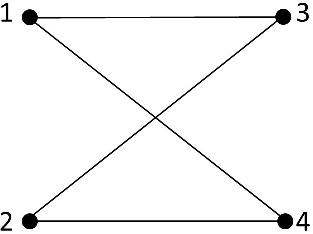}　\quad
\includegraphics[width = 3 cm]{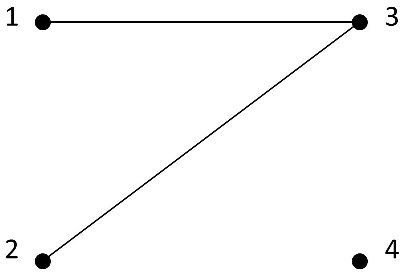}

Figure 1.3\quad Graphs $\Gamma(R_1), \Gamma(R_2)$
\end{center} 
\end{ex*}

The following lemma can easily be proved, see \cite[Section 6]{em}. 

\begin{lemm}\label{lem3}
For $R \subset \Phi$ with $R \cap (-R) = \emptyset$, the graph $\Gamma(R)$ has no cycle if and only if elements of $R$ are linearly independent, and $\Gamma(R)$ is connected if and only if the linear subspace spanned by the elements of $R$ is of codimension one in $\mathbb{R}^n$. Therefore, $\Gamma(R)$ is tree if and only if $R$ has $n-1$ elements and they are linearly independent. 
\end{lemm}

Henceforth our matroid polytope $\DM$ in $\mathbb{R}^n$ is assumed to be of maximal dimension, i.e. $\dim \DM=n-1$ unless otherwise stated. 
Let $v$ be a vertex of $\DM$. By Proposition~\ref{prop:1-1}, $v=\delta_J$ for some $J\in \mathcal{M}$ and any edge vector emanating from $v$ is of the form $\delta_I-\delta_J$ for $I\in\mathcal{M}$ with $|I\cap J|=k-1$. Note that $\delta_I-\delta_J=\delta_i-\delta_j$ where $\{i\}=I\backslash J$ and $\{j\}=J\backslash I$. Based on this observation, for $v=\delta_J\in V(\DM)$ we define 
\begin{equation} \label{eq:1-2}
\Rv(\DM)=\{ \delta_i-\delta_j\mid \{i\}=I\backslash J,\ \{j\}=J\backslash I,\ \delta_I\delta_J\in E(\DM)\}.
\end{equation}

\begin{lemm}\label{lem6}
The graph $\Gamma(\Rv(\DM))$ is a connected bipartite graph with parts $J$ and $[n] \setminus J$. 
It is a tree if and only if $\DM$ is simple at $v \in V(\Delta_\mathcal{M})$. 
\end{lemm}

\begin{Proof}
It follows from \eqref{eq:1-2} and the definition of $\Gamma(R)$ that any edge of $\Gamma(\Rv(\DM))$ has one endpoint in $J$ and the other endpoint in $[n]\backslash J$. This together with Lemma~\ref{lem3} implies the lemma. 
\QED
\end{Proof}

For a vertex $v=\delta_J\in \Delta_{n,k}$ and a positive integer $s$ we introduce two notations: 
\begin{equation}
\begin{split}
\Wvs &= \{L \subset [n] \mid |L \cap J| = |L \setminus J| = s\},\\
\Vvs&=\{\delta_K \in V(\Delta_{n,k}) \mid K\subset [n],\ |K|=k,\ |J \cap K| = k-s\}.
\end{split}
\notag
\end{equation}
There is a bijection 
\begin{equation} \label{eq:1-3}
\varphi\colon \Wvs \to \Vvs.
\end{equation}
Indeed, to $L\in \Wvs$ we associate the symmetric difference of $L$ and $J$ as $K$, that is, $K=(J\backslash L)\cup (L\backslash J)$, and this induces the bijection $\varphi$. 
Note that $\Vvone$ is the set of vertices of $\Delta_{n,k}$ that are connected to $v=\delta_J$ by an edge of $\Delta_{n,k}$. If $i \notin J$ and $j \in J$, then $\{i, j\} \in \Wvone$ and 
\[
\varphi(\{i,j\})=\delta_J +(\delta_i- \delta_j)\in\Vvone\subset V(\Delta_{n,k}).
\]
With this understood, we have 

\begin{lemm}\label{lem7}
$\varphi(\{i,j\})$ is a vertex of $\DM$ if and only if $i$ and $j$ are joined by an edge of $\Gamma(\Rv(\Delta_\mathcal{M}))$. 
\end{lemm}

\begin{Proof}
The lemma follows from the following observation: 
\begin{align}
\delta_J+(\delta_i - \delta_j)\in V(\DM) &\Leftrightarrow \delta_i - \delta_j \in \Rv(\DM) \notag \\
&\Leftrightarrow \mbox{$i$ and $j$ are joined by an edge of } \Gamma(\Rv(\DM)). \notag\QED
\end{align}
\end{Proof}

\begin{ex*}
Let $S = [4]$, $k = 2$, $J = \{1, 2\}$ and consider the following matroid polytope $\Delta$:

\begin{center}
\includegraphics[width = 7 cm]{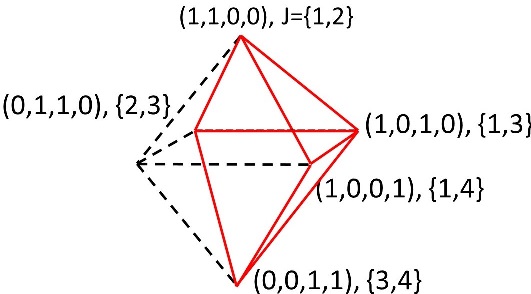}

Figure 1.4\quad An example of a matroid polytope with the vertex $\delta_J$
\end{center}
Then, for $v = \delta_J = (1,1,0,0)$, the graph $\Gamma = \Gamma(\Rv(\Delta))$ associated to $\Delta$ is the following: 
\begin{center}
\includegraphics[width = 3 cm]{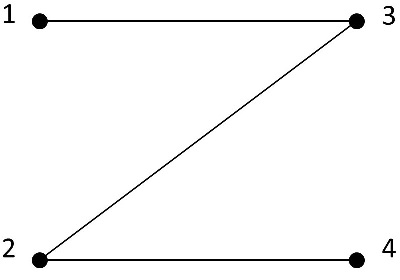}

Figure 1.5\quad The bipartite graph $\Gamma$ associated to $\Delta$
\end{center}
Indeed, since 
\begin{align}
(1,0,0,1)-(1,1,0,0) = (0,-1,0,1) &= \delta_4 - \delta_2 \notag \\
(1,0,1,0)-(1,1,0,0) = (0,-1,1,0) &= \delta_3 - \delta_2 \notag \\
(0,1,1,0)-(1,1,0,0) = (-1,0,1,0) &= \delta_3 - \delta_1, \notag 
\end{align}
the vertices $4$ and $2$, $3$ and $2$, $3$ and $1$ are respectively joined by an edge. 
\end{ex*}

We recall that if $\Gamma$ is a graph and $U$ is a subset of the vertex set of $\Gamma$, then the induced subgraph, denoted $\Gamma|U$, is the subgraph of $\Gamma$ whose vertex set is $U$ and whose edge set consists of all of the edges in $\Gamma$ that have both endpoints in $U$. The following lemma plays a key role to prove Theorem $\ref{theo1}$. 

\begin{lemm}\label{lem8} 
Let $v = \delta_J \in V(\DM)$ as before and $U \in \Wvtwo$. If the induced subgraph $\Gamma(\Rv(\DM))|U$ is isomorphic to one of the following two graphs:
\begin{center}
\includegraphics[width = 5 cm]{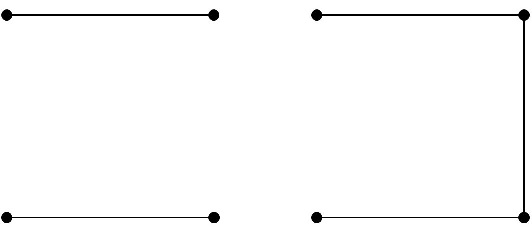}

Figure 1.6 \quad The induced subgraphs $\Gamma(\Rv(\DM))|U$
\end{center}
then $\varphi(U)$ is a vertex of $\DM$, where $\varphi$ is the bijection in \eqref{eq:1-3} 
\end{lemm}

\begin{Proof}
We may assume $U = \{1,2,3,4\}$ and $J = \{1,2,5,6,\dots, k+2\}$ through a permutation on $[n]$ if necessary. Let $J_0 = \{5,6,\dots, k+2\}$.
Consider the half space
\begin{equation}
\{ (\xi_1,\dots,\xi_n)\in \mathbb{R}^n\mid \xi_5 + \xi_6 + \cdots + \xi_{k+2} - \xi_{k+3} - \cdots - \xi_{n} \leq k-2\} \notag
\end{equation}
and let $H$ be the hyperplane in $\mathbb{R}^n$ that is the boundary of the above half space. Since the half space above contains $\Delta_{n,k}$, the intersection $\Delta_{n,k}\cap H$ is a face of $\Delta_{n,k}$, denoted by $f$. Note that 
\begin{equation} \label{eq:1-4}
\text{$f$ is the convex hull of $\{\delta_A+\delta_{J_0} \mid A \subset U=\{1,2,3,4\},\ |A| = 2\}$,} 
\end{equation}
so $f$ is isomorphic to $\Delta_{4,2}$. Similarly, since $\DM$ is contained in $\Delta_{n,k}$, the intersection $\DM\cap H$ is a face of $\DM$. Here 
\[
\DM \cap f=\DM\cap (\Delta_{n,k}\cap H)=(\DM\cap \Delta_{n,k})\cap H=\DM\cap H,
\] 
so $\DM\cap f$ is a face of $\DM$ and hence we have 
\[
V(\DM \cap f) \subset V(\DM) \subset V(\Delta_{n,k}), \quad 
E(\DM \cap f) \subset E(\DM) \subset E(\Delta_{n,k}). 
\]
Therefore, 
\begin{equation} \label{eq:1-5}
\text{$\DM\cap f$ is a matroid polytope}
\end{equation}
by Proposition~\ref{prop:1-1}. 

The $v = \delta_{J}$ is a vertex of $\DM$ by assumption and also a vertex of $f$ by \eqref{eq:1-4} Therefore it is a vertex of $\DM\cap f$ since $\DM \cap f$ is a face of $\DM$. 
Since $U=\{1,2,3,4\}$ and $J=\{1,2,5,6,\dots,k+2\}$ as remarked at the beginning of the proof, we have 
\[
\varphi(U)=\delta_{\{3,4,5,6,\dots,k+2\}}=\delta_{\{3,4\}}+\delta_{J_0}
\]
and it is a vertex of $f$ by \eqref{eq:1-4}. We shall show that $\delta_{\{3,4\}}+\delta_{J_0}$ is a vertex of $\DM\cap f$ under the assumption in the lemma.

\begin{center}
\includegraphics[width=5cm]{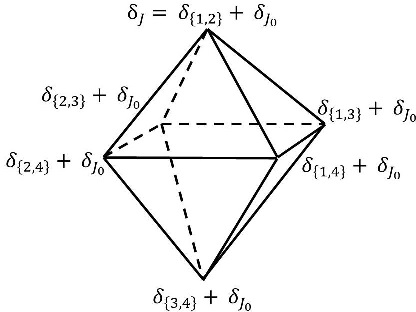}

Figure 1.7 \quad The face $f=\Delta_{4,2}$
\end{center}

{\bf Case 1.} 
The case where the induced subgragh $\Gamma(\Rv(\DM))|U$ is isomorphic to the graph on the left hand side of Figure 1.6. In this case we may assume without loss of generality that $\{1, 3\}$ and $\{2, 4\}$ are edges of $\Gamma(\Rv(\DM))$ while $\{1,4\}$ and $\{2,3\}$ are not. Therefore it follows from Lemma~\ref{lem7} that 
\[
\delta_J + \delta_3 - \delta_1 = \delta_{\{2,3\}} + \delta_{J_0} \quad\text{and}\quad \delta_J + \delta_4 - \delta_2 = \delta_{\{1,4\}} + \delta_{J_0}
\]
are vertices of $\DM$ while 
\[
\delta_J + \delta_4 - \delta_1 = \delta_{\{2,4\}} + \delta_{J_0}\quad\text{and}\quad \delta_J + \delta_3 - \delta_2 = \delta_{\{1,3\}} + \delta_{J_0}
\]
are not. Since the four elements above are vertices of $f$ by \eqref{eq:1-4}, the former two elements above are vertices of $\DM\cap f$ while the latter two elements above are not. 
Then $\delta_{\{3,4\}} + \delta_{J_0}$ must be a vertex of $\DM \cap f$ by \eqref{eq:1-5}, see the figure on the left side of Figure 1.8. 

\begin{center}
$\delta_J$ \hspace{1in} $\delta_J$

\includegraphics[width=6cm]{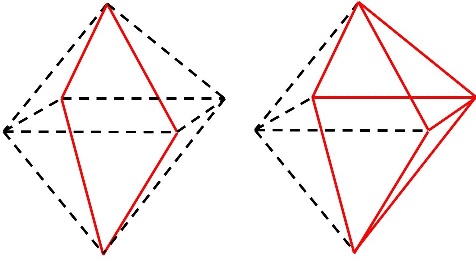}

\hspace{-0.6in} Figure 1.8\quad Case $1$ \hspace{0.5in} Case $2$
\end{center}

{\bf Case 2.} 
The case where the induced subgraph $\Gamma(\Rv(\DM))|U$ is isomorphic to the graph on the right side of Figure 1.6. In this case we may assume without loss of generality that $\{1, 3\}$, $\{2, 4\}$ and $\{1,4\}$ are edges of $\Gamma(\Rv(\DM))$ while $\{2,3\}$ is not. Then the same argument as in Case 1 shows that $\delta_{\{2,3\}} + \delta_{J_0}$, $\delta_{\{1,4\}} + \delta_{J_0}$ and $\delta_{\{2,4\}} + \delta_{J_0}$ are vertices of $\DM\cap f$ while $\delta_{\{1,3\}}+\delta_{J_0}$ is not. Then $\delta_{\{3,4\}} + \delta_{J_0}$ must be a vertex of $\DM \cap f$ by \eqref{eq:1-5} as before, see the figure on the right side of Figure 1.8. \QED
\end{Proof}

Since it suffices to prove Theorem~\ref{theo1} when the matroid polytope $\DM$ is of maximal dimension by Proposition~\ref{prop:1-2}, the theorem follows from the following. 

\begin{theo}\label{thm2}
Let $k, n-k \geq 2$ and $\Delta$ be a matroid polytope of a matroid of rank $k$ on $[n]$ such that $\dim\Delta=n-1$. If $v$ is a simple vertex of $\Delta$, then there is a non-simple vertex of $\Delta$ adjacent to $v$. 
\end{theo}

\begin{Proof}
Let $\Gamma$ denote $\Gamma(\Rv(\Delta))$ and $v=v_J$. 
Since $\Delta$ is of maximal dimension and $v=v_J$ is a simple vertex, $\Gamma$ is a tree with vertices $[n]$ and bipartite with parts $J$ and $[n]\backslash J$ by Lemma~\ref{lem6}. Both $|J|=k$ and $|[n]\backslash J|=n-k$ are more than one by assumption, so $\Gamma$ has an edge, say $\{i,j\}$, whose endpoints $i$ and $j$ are not endpoints of the graph $\Gamma$. We may assume $i \in [n] \setminus J$ and $j \in J$ since $\Gamma$ is a bipartite graph with parts $J$ and $[n]\backslash J$. Note that $\{i,j\}\in \Wvone$ (see the paragraph above (5)). 

Let $v'=\varphi(\{i,j\})\in \Vvone$ where $\varphi$ is the bijection in \eqref{eq:1-3}. The $v'$ is a vertex of $\Delta$ by Lemma \ref{lem7} and adjacent to $v$ since $v'\in \Vvone$. 
We shall prove that $v'$ is the desired non-simple vertex of $\Delta$. Since $v'=\varphi(\{i,j\})$, we have $v'=\delta_J+\delta_i-\delta_j$. Let $J'=(J \setminus \{j\}) \cup \{i\}$. Then $v'=\delta_{J'}$. Let $V'$ denote the set of vertices of $\Delta$ adjacent to $v'=\delta_{J'}$. Since $\dim\Delta=n-1$, we will prove that $|V'|\ge n$ which means that $\Delta$ is not simple at $v'$. 

Note that 
\[
V'\subset V_{v',1}=\{\delta_K \mid K \subset [n],\ |J' \cap K|=k-1\} \subset V(\Delta_{n,k})
\]
and $V_{v',1}$ is the set of vertices of $\Delta_{n,k}$ adjacent to $v'$. One can write $K=(J' \setminus \{p\}) \cup \{q\}$ with some $p \in J'$ and $q \in [n] \setminus J'$. We take four cases: 
\[
\text{$(a) \ p=i, q=j, \ \ (b) \ p=i, q\not=j, \ \ (c) \ p\not=i, q=j, \ \ (d) \ p\not=i, q\not=j $.}
\] 
In each case, $K$ is written as follows:
\begin{equation}
K=
\begin{cases}
(a) \ J \\
(b) \ (J\backslash \{j\}) \cup \{q\} & \ (q \in ( [n]\backslash J)\backslash \{i\}) \\
(c) \ (J\backslash \{p\}) \cup \{i\} & \ (p \in J\backslash \{j\}) \\
(d) \ (J\backslash L_1) \cup L_2 & \ (j \in L_1 \subset J,\ i \in L_2 \subset [n]\backslash J,\ |L_1|=|L_2|=2).
\end{cases} \notag
\end{equation}
In case $(a)$ we have $\delta_K=v$. In case $(b)$ or $(c)$ we have $|J \cap K|=k-1$ and in case $(d)$ we have $|J \cap K|=k-2$. Therefore, through the symmetric difference $(K\backslash J)\cup (J\backslash K)$, $K$ corresponds to the following element in case $(b)$, $(c)$ or $(d)$:
\begin{equation}
\begin{cases}
(b) \ L \in \Wvone \quad (L \cap \{i, j\} =\{j\}) \\
(c) \ L \in \Wvone \quad (L \cap \{i, j\} =\{i\}) \\
(d) \ U \in \Wvtwo \quad (\{i, j\} \subset U ).
\end{cases} \notag
\end{equation} 
Note that $\delta_K$ is $\varphi(L)$ in case $(b)$ or $(c)$ and $\varphi(U)$ in case $(d)$. 

We shall observe which $\delta_K$ is a vertex of the polytope $\Delta$. In case $(a)$, $\delta_K=v$ that is a vertex of $\Delta$. In case $(b)$ or $(c)$, it follows from Lemma~\ref{lem7} that $\delta_K$ is a vertex of $\Delta$ if $L$ is an edge of the graph $\Gamma$. Note that such an edge $L$ meets the edge $\{i,j\}$ at the vertex $j$ or $i$ since $L\cap\{i,j\}=\{j\}$ or $\{i\}$. 
In case $(d)$, write $U$ as $\{i,j,x,y\}$ and look at the induced subgraph $\Gamma|U$. If $\{x,y\}$ is an edge of $\Gamma$ not adjacent to $\{i, j\}$, then $\Gamma|U$ is isomorphic to the graph on the left hand side of Figure 1.6 because $\Gamma$ is a tree. Therefore $\delta_K$ is a vertex of $\Delta$ by Lemma \ref{lem8} in this case. 
Thus each edge of $\Gamma$ produces a vertex of $\Delta$ adjacent to the vertex $v'$ and these vertices are all distinct. Since $\Gamma$ is a tree with vertices $[n]$, $\Gamma$ has $n-1$ edges; so $|V'|\ge n-1$. 

\begin{center}
\includegraphics[width=5.5cm]{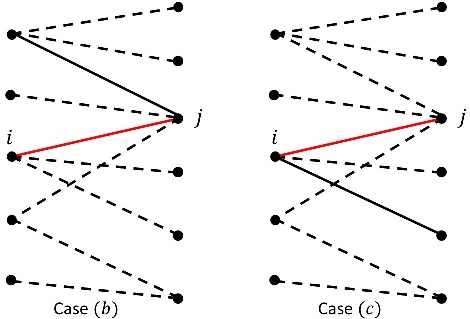} \quad
\includegraphics[width=5.5cm]{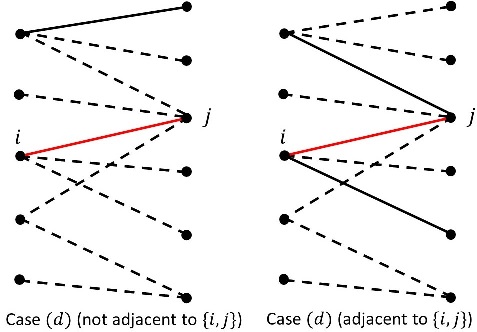}

Figure 1.9 \quad Examples of edges of $\Gamma$ in case $(b), (c)$ or $(d)$
\end{center}

We shall prove that there is another element in $V'$. 
Since the endpoints $i$ and $j$ of the edge $\{i,j\}$ are not the endpoints of $\Gamma$ by our choice of the edge $\{i,j\}$, there are edges $\{i, a\}$ and $\{j, b\}$ of $\Gamma$ with some $a,b\in [n]$. 
More precisely, $a\in J$ and $b\in [n]\backslash J$ because $\Gamma$ is a bipartite graph with parts $J$ and $[n]\backslash J$ and $i\in [n]\backslash$ and $j\in J$. Therefore $\{i, j, a, b\} \in\Wvtwo$. We take $U=\{i, j, a, b\}$. Since $\Gamma|U$ has at least three edges $\{i,j\}, \{i,a\}, \{j,b\}$ and $\Gamma$ is tree, $\Gamma|U$ must be isomorphic to the graph on the right side of Figure 1.6. Therefore, $\delta_K=\varphi(U)$ is a vertex of $\Delta$ by Lemma~\ref{lem8} that is the desired another element of $V'$. \QED
\end{Proof}

\section{Smooth torus orbit closures in the Grassmannians}

For a subset $J \subset [n]$, let $\mathbb{C}^J$ be the set of complex-valued functions on $J$.
We identify $\mathbb{C}^{[n]}$ with $\mathbb{C}^n$. 
Let $G_k(E)$ denote the Grassmannian consisting of $k$-dimensional linear subspaces of a complex vector space $E$. 
When $(J_1,J_2,\dots,J_r)$ is a partition of $[n]$, we naturally identify $\mathbb{C}^{J_1} \times \mathbb{C}^{J_2} \times \cdots \times \mathbb{C}^J_r$ with $\mathbb{C}^n$ and regard $G_{k_1}(\mathbb{C}^{J_1}) \times G_{k_2}(\mathbb{C}^{J_2}) \times \cdots \times G_{k_r}(\mathbb{C}^{J_r})$ as a subvariety of $G_k(\mathbb{C}^n)$ accordingly, where $(k_1,k_2,\dots,k_r)$ is a partition of the integer $k$ with $k_1 \leq |J_1|,k_2 \leq |J_2|, \ldots , k_r \leq |J_r|$. 

The general linear group ${\rm GL}(\mathbb{C}^n)$ naturally acts on $G_k(\mathbb{C}^n)$. Let $T$ be the subgroup of ${\rm GL}(\mathbb{C}^n)$ consisting of diagonal non-singular matrices. We denote the $T$-orbit of $X$ by $TX$ and its closure in $G_k(\mathbb{C}^n)$ by $\overline{TX}$. It is known that $\overline{TX}$ is a toric variety (see \cite[p.135, 1.1]{gs} for example). Our main result in this section is the following. 

\begin{theo}\label{main2}
Let $X \in G_k(\mathbb{C}^n)$. If the torus orbit closure $\overline{TX}$ is smooth, then there exist a partition $(J_1,J_2,\dots,J_r)$ of the set $[n]$ and a partition $(k_1,k_2,\dots,k_r)$ of the integer $k$ with $k_i \in \{0, 1, |J_i|-1, |J_i|\}$ such that
\begin{equation}
\overline{TX}=G_{k_1}(\mathbb{C}^{J_1}) \times G_{k_2}(\mathbb{C}^{J_2}) \times \cdots \times G_{k_r}(\mathbb{C}^{J_r}), 
\end{equation}
in particular, $\overline{TX}$ is a product of complex projective spaces. 
\end{theo}

We recall a moment map on $G_k(\mathbb{C}^n)$ which reduces Theorem~\ref{main2} to Theorem \ref{theo1}. 
Let $I_{n,k}$ denote the set of $k$-element subsets of $[n]$. Given $X \in G_k(\mathbb{C}^n)$, we represent it as an $n \times k$ matrix, say $A$, and for any $J=\{j_1,j_2,\ldots, j_k\}\in I_{n,k}$, define $P^J(X)$ to be the minor of the matrix consisiting of the rows of $A$ with indices $j_1, j_2, \ldots, j_k$. 
The moment map $\mu\colon G_k(\mathbb{C}^n)\to \mathbb{R}^n$ is defined to be 
\begin{equation}
\mu(X) = \frac{\sum_{J \in I_{n,k}}|P^J(X)|^2\cdot \delta_J}{\sum_{J \in I_{n,k}}|P^J(X)|^2}.\notag
\end{equation}
The image $\mu(G_k(\mathbb{C}^n))$ is the hypersimplex $\Delta_{n,k}$, that is the convex hull of $\{ \delta_J \mid J \in I_{n,k} \}$ in $\mathbb{R}^n$. More generally, it is known that $\mu(\overline{TX})$ is a convex polytope in $\mathbb{R}^n$ with vertex set $\{\delta_J\mid P^J(X) \neq 0\}$ and the map $\mu \colon \overline{TX}\to \mu(\overline{TX})$ induces a one-to-one correspondence between the $p$-dimensional orbits of $T$ in $\overline{TX}$ and the open $p$-dimensional faces of $\mu(\overline{TX})$ (see \cite[Theorem 1 in p.139]{gs} for example). 

The key fact on the moment map used in our argument is the following (see \cite[Proposition 1 in p.150]{gs} for example).

\begin{prop} 
\label{convexity}
Let $X\in G_k(\mathbb{C}^n)$. Then $\overline{TX}$ is smooth if and only if $\mu(\overline{TX})$ is a simple polytope.
\end{prop}

\begin{re*}  The theory of toric varieties says that if the moment polytope $\mu(\overline{TX})$ is simple, then $\overline{TX}$ is an orbifold.   In our case, the primitive edge vectors of $\mu(\overline{TX})$ are of the form $\pm(\delta_i-\delta_j)$, so simpleness of $\mu(\overline{TX})$ implies that $\overline{TX}$ is non-singular (see \cite[Lemma 6.6]{em} for example).    
\end{re*}

When $(n,k)=(n, 0)$ or $(n,n)$, $G_k(\mathbb{C}^n)$ is a point and when $(n,k)=(n,1)$ or $(n, n-1)$, $G_k(\mathbb{C}^n)$ is the complex projective space of complex dimension $n-1$. Therefore, the following lemma is obvious. 

\begin{lemm}\label{maxorbit-maximage}
Let $X \in G_k(\mathbb{C}^n)$.  When $(n,k)=(n, 0), (n,1), (n, n-1)$ or $(n,n)$, we have 
\begin{center}
$\overline{T X}=G_k(\mathbb{C}^n) \Leftrightarrow \mu(\overline{T X})=\Delta_{n,k}$.
\end{center}
\end{lemm}

Suppose that $(J_1,J_2,\dots,J_r)$ is a partition of the set $[n]$ and $(k_1,k_2,\dots,k_r)$ is a partition of the integer $k$ with $k_1 \leq |J_1|,k_2 \leq |J_2|, \ldots , k_r \leq |J_r|$. 
Let $T_i$ be the subgroup of ${\rm GL}(\mathbb{C}^{J_i})$ consisting of non-singular diagonal matrices. 
We naturally identify $T_1\times T_2\times \cdots \times T_r$ with $T$. Then for 
$X=(X_1,X_2,\dots,X_r)\in G_{k_1}(\mathbb{C}^{J_1}) \times G_{k_2}(\mathbb{C}^{J_2}) \times \cdots \times G_{k_r}(\mathbb{C}^{J_r})$, one can easily see that 
\[
\overline{TX}=\overline{T_1 X_1} \times \overline{T_2 X_2 } \times \cdots \times \overline{T_r X_r }\quad\text{and}\quad 
\mu(\overline{T X})=\mu_1(\overline{T_1 X_1})
\times \mu_2(\overline{T_2 X_2}) \times \cdots \times \mu_r(\overline{T_r X_r })
\]
where $\mu_i$ denotes the moment map on $G_{k_i}(\mathbb{C}^{J_i})$ for each $i=1,2,\dots,r$. 
Conversely, the following holds. 

\begin{lemm}\label{split}
Let $X\in G_k(\mathbb{C}^n)$. If 
\begin{center}
$\mu(\overline{T X})=\Delta_1 \times \Delta_2 \times \cdots \times \Delta_r\subset \mathbb{R}^{J_1}\times \mathbb{R}^{J_2}\times \dots \times \mathbb{R}^{J_r}=\mathbb{R}^n$
\end{center}
where $\Delta_i$ is a matroid polytope of some matroid of rank $k_i$ on $J_i$, 
then $X \in G_k(\mathbb{C}^n)$ belongs to $G_{k_1}(\mathbb{C}^{J_1}) \times G_{k_2}(\mathbb{C}^{J_2}) \times \cdots \times G_{k_r}(\mathbb{C}^{J_r})$ and if we write $X$ as $(X_1, X_2, \ldots, X_r)$ with $X_i\in G_{k_i}(\mathbb{C}^{J_i})$ for each $1 \leq i \leq r$, then $\mu_i(\overline{T_i X_i})=\Delta_i$.
\end{lemm}

We will give the proof of this lemma later and complete the proof of our main theorem. 

\begin{Proof}[Proof of Theorem \ref{main2}] 
Suppose that $\overline{T X}$ is smooth. Then $\mu(\overline{T X})$ is a simple polytope by Proposition~\ref{convexity}. Therefore, it follows from Theorem \ref{theo1} that there exist a partition $(J_1,J_2,\dots,J_r)$ of the set $[n]$ and a partition $(k_1,k_2,\dots,k_r)$ of the integer $k$ with $k_i \in \{0,1,|J_i|-1,|J_i|\}$ such that 
\begin{equation}
\mu(\overline{T X}) = \Delta_{J_1,k_1} \times \Delta_{J_2,k_2} \times \cdots \times \Delta_{J_r,k_r}. \notag
\end{equation}
From Lemma \ref{split}, $X=(X_1, X_2, \ldots, X_r) \in G_{k_1}(\mathbb{C}^{J_1}) \times G_{k_2}(\mathbb{C}^{J_2}) \times \cdots \times G_{k_r}(\mathbb{C}^{J_r})$ 
and $\mu_i(\overline{T_i X_i})=\Delta_{J_i, k_i}$ for each $1 \leq i \leq r$. Since $k_i \in \{0, 1, |J_i|-1, |J_i|\}$, we have $\overline{T_iX_i}=G_{k_i}(\mathbb{C}^{J_i})$ for each $1\le i\le r$ by Lemma \ref{maxorbit-maximage}. Thus we have 
\begin{center}
$\overline{T X}=\overline{T_1X_1}\times \overline{T_2X_2}\times \cdots \times \overline{T_rX_r}=G_{k_1}(\mathbb{C}^{J_1}) \times G_{k_2}(\mathbb{C}^{J_2}) \times \cdots \times G_{k_r}(\mathbb{C}^{J_r})$,
\end{center}
proving the theorem. \QED
\end{Proof}

\vspace{0.3in}

We shall prove Lemma \ref{split}.

\begin{Proof}[Proof of Lemma \ref{split}]
It suffices to prove Lemma \ref{split} for $r=2$. By permuting coordinates of $\mathbb{C}^n$, we may assume that $J_1=[|J_1|]$ and $J_2=[n] \setminus [|J_1|]$. Let $X\in G_k(\mathbb{C}^n)$. Suppose that 
\begin{equation}\label{assumption}
\mu(\overline{T X})=\Delta_1 \times \Delta_2 \subset \mathbb{R}^{J_1} \times \mathbb{R}^{J_2}.
\end{equation}
If $X$ is represented by a $n \times k$ matrix of the form
\begin{equation}
\left(
\begin{array}{cc}
A_1 & 0 \\
0 & A_2 \\
\end{array}
\right), \notag
\end{equation}
where each $A_i$ is a matrix of type $|J_i| \times k_i$, then we can take $X_i \in G_{k_i}(\mathbb{C}^{J_i})$ represented by the matrix $A_i$ and we have $\mu_i(\overline{T_i X_i})=\Delta_i$ for each $i$. Let $X$ be represented by a $n \times k$ matrix of rank $k$
\begin{equation}
A=
\left(
\begin{array}{c} 
\alpha_1\\ \alpha_2 \\ \vdots \\ \alpha_n 
\end{array} 
\right), \notag
\end{equation} 
where $\alpha_1, \alpha_2 \ldots \alpha_n$ are row vectors. Let $V$ be the vector space generated by $\alpha_1, \alpha_2, \ldots, \alpha_n$, $V_1$ be the subspace generated by $\alpha_1, \alpha_2, \ldots, \alpha_{|J_1|}$ and $V_2$ the subspace generated by $\alpha_{|J_1|+1}, \alpha_{k_1+2}, \ldots, \alpha_n$. We shall prove that $V=V_1 \oplus V_2$. If this claim is proved, then there is a non-singular $k \times k$ matrix $g$ such that 
\begin{equation}
Ag=
\left(
\begin{array}{cc}
A_1 & 0 \\
0 & A_2 \\
\end{array}
\right), \notag
\end{equation}
which completes the proof since the matrix $A g$ also represents $X$. Let a matroid $\mathcal{M}_i$ on $J_i$ of rank $k_i$ correspond to the matroid polytope $\Delta_i$ for each $i$. Then assumption \eqref{assumption} means that 
\begin{equation}\label{eqfor-2-4}
\{ B \in I_{n, k} | \ P^B(X) \neq 0 \}=\{B_1 \sqcup B_2 \ | \ B_1 \in \mathcal{M}_1, \\ B_2 \in \mathcal{M}_2 \}. 
\end{equation}
Note that $B_i \subset J_i$ for each $i$. Since the $n \times k$ matrix $A$ is of rank $k$, $V$ equals the vector space of row vectors with $k$ components. From \eqref{eqfor-2-4}, any basis of $V$ consisting of $k$ row vectors of $A$ consists of $k_1$ row vectors of $A$ with indices in $J_1$ and $k_2$ row vectors of $A$ with indices in $J_2$. Thus $k_1 \leq {\rm dim} \ V_1$ and $k_2 \leq {\rm dim} \ V_2$. Suppose $k_1 < {\rm dim} \ V_1$ (resp. $k_2 < {\rm dim} \ V_2$). Then there is a basis of $V$ consisting of $k$ row vectors of $A$ such that some vectors of the basis form a basis of $V_1$ (resp. $V_2$). However, for any basis of $V$ consisting of $k$ row vectors of $A$, the number of members of the basis with indices in $J_1$ (resp. $J_2$) must be $k_1$ (resp. $k_2$). This is a contradiction. Thus $k_1 = {\rm dim} \ V_1$ and $k_2 = {\rm dim} \ V_2$. \QED
\end{Proof}

\textsc{Division of Mathematics \& Physics, Graduate School of Science, Osaka City University, 3-3-138 Sugimoto, Sumiyoshi-ku, Osaka 558-8585, Japan}

\textit{E-mail address} : \texttt{mathlibrary0824@gmail.com} 

\hspace{1.05in}\texttt{m18sa005@du.osaka-cu.ac.jp}

\end{document}